\documentclass{amsart}
\usepackage{amsfonts}
\begin{document}

\newtheorem{Theorem}{Theorem}[section]
\newtheorem{Example}{Example}[section]
\newtheorem{Remark}{Remark}[Theorem]
\newtheorem{Lemma}{Lemma}[Theorem]
\newtheorem{Definition}{Definition}[section]

\newcommand{\Imag}{\operatorname{Im}}
\newcommand{\Id}{\operatorname{Id}}
\newcommand{\End}{\operatorname{End}}
\newcommand{\Tr}{\operatorname{Tr}}

\def\qedbox{\hbox{$\rlap{$\sqcap$}\sqcup$}}
\def\qed{\nobreak\hfill\penalty250 \hbox{}\nobreak\hfill\qedbox}
\def\e{\varepsilon}
\def\diag{\mbox{diag}}
\def\dim{\mbox{dim}}
\def\arg{\mbox{arg}}
\def\tr{\mbox{tr}}
\def\e{\epsilon}
\def\ve{\varepsilon}
\def\RR{\mathbb{R}}
\def\HH{\mathbb{H}}
\def\FF{\mathbb{F}}
\def\NN{\mathbb{N}}
\def\ZZ{\mathbb{Z}}
\def\II{\mathbb{1}}

\def\PHPn{\tilde{\HH} P^n}
\def\PHn{{\tilde{\HH}}^n}
\def\CC2n{{\CC}^{2n}}
\def\PH{\tilde{\HH}}
\def\CC{\mathbb{C}}
\def\H1{{\tilde{\HH}}_1}

\def\cal#1{{\mathfrak #1}}

\title{Para-quaternionic reduction}
\author{S. Vukmirovi\' c}
\address{Faculty of mathematics, University of Belgrade, Studenski trg 16, p.p.
550, 11 000 Belgrade, Yugoslavia}
\email{vsrdjan@matf.bg.ac.yu }
\thanks{Work completed during the tenure of a University of Hull graduate teaching assistantship of second author}
\date{}
\begin{abstract}
The pseudo-Riemannian manifold $M=(M^{4n},g), n \geq 2$ is
para-qua\-te\-rni\-onic K\" ahler if $hol(M) \subset sp(n, \RR)
\oplus sp(1, \RR ).$ If $hol(M) \subset sp(n, \RR),$ than the
manifold $M$ is called para-hyperK\" ahler.  The other possible
definitions of these manifolds use certain parallel
para-quaternionic structures in $\End (TM),$ similarly to the
quaternionic case. In order to relate these different definitions
we study para-quaternionic algebras in details. We describe the
reduction method for the para-quaternionic K\" ahler and
para-hyperK\" ahler manifolds  and give some  examples. The
decomposition of a curvature tensor of the para-quaternionic type
is also described.
\end{abstract}
\maketitle

\section{Introduction}
\noindent In the paper we try to develope  theory of the
para-quaternionic structures on a manifold. Algebra of
para-quaternions is  known but rarely used so far. A  nice
overview of para-quaternions and related classical geometries is
given in \cite{R} (para-quaternions were referred as generalized
quaternions).
 The  para-quaternionic sectional
curvature of the para-quaternionic projective space $\PHPn$ has
been studied in \cite{B}. The para-hyperK\" ahler manifolds have
been  studied in \cite{H1, Kam1, Kam2}, where they are referred as
neutral hyperK\" ahler manifolds. Quite recently, the author have
learnt that, independently, the notion of a para-quaternionic K\"
ahler manifold has been developed in \cite{GMV}. Although the
basic definitions and conclusions are the same, the investigations
in \cite{GMV} and this paper go in different directions. In
\cite{GMV}  the para-quaternionic sectional curvature of the
para-quaternionic K\" ahler manifold and its relations to the
Osserman condition have been studied.

\noindent In this paper we systematically study the para-hyperK\"
ahler and para-quaternionic K\" ahler manifolds. We characterize
them as  the manifolds with holonomies contained in $sp(n, \RR )
\oplus sp(1, \RR)$ and  $sp(n, \RR) ,$ respectively. They are
necessarily of the dimension $4n, n \geq 2$ and of the signature
$(2n, 2n).$ Since $sp(n, \RR)$ are $sp(n)$ are  real forms of the
same complex Lie algebra,  para-hyperK\" ahler and
para-quaternionic K\" ahler manifolds enjoy similar properties to
the quaternionic counterparts:   hyperK\" ahler and quaternionic
K\" ahler manifolds. Some  facts concerning the quaternionic
geometry given in \cite{Sal}  can be carried in the
para-quaternionic setting by means of the complexification. The
classification of para-hyperK\" ahler symmetric spaces is given in
\cite{ABCV}. The decomposition of the space of curvature tensors
of a para-quaternionic type, given here,  is analogous to the
decomposition   in the quaternionic case (see \cite{AM}). The
reduction methods for the hyperK\" ahler and quaternionic K\"
ahler manifolds are well known (see \cite{HKLR, GL}). Here we
describe the reduction methods for the para-hyperK\" ahler and
para-quaternionic K\" ahler manifolds.

\noindent However, many questions  remained open in the
para-quaternionic case. It is known (see \cite{BF, K}) that the
cotangent bundle $T^*M$ of a K\" ahler manifold $M$ is at least
locally a hyperK\" ahler manifold. It is an interesting question
what the natural space, on which para-hyperK\" ahler structures
can arise, is. In \cite{AS} quaternionic K\" ahler manifolds have
been studied as certain quotients of special types of the
hyper-K\" ahler manifolds. Does  similar relation between
para-hyperK\" ahler and para-quaternionic K\" ahler manifolds
exist?

 \noindent The structure of the paper is as follows. In
Section \ref{linear} we establish a basic notation and  give basic
properties of the para-quaternions and Lie algebras related to
them. Subsection \ref{grassman} enables us to use the Grassman
formalism in the para-quaternionic vector space. In Section
\ref{manifolds} we define the notions of the para-hyperK\" ahler
and para-quaternionic K\" ahler manifold and give their
characterizations in terms of the algebra of holonomy (Theorems
\ref{charaPHK} and \ref{charaPQK}). A detailed analysis of the
para-quaternionic projective space is given in Subsection
\ref{projspace}. In Section \ref{decomposition} we give the
decompositions of spaces of the curvature tensors of the type
$gl_{2n}( \RR ) \oplus sp(1, \RR)$ (Theorem \ref{decompGL}) and
$sp(n, \RR ) \oplus sp(1, \RR)$ (Theorem \ref{decompSP}). Finally,
in Section \ref{reduction} we describe the reduction techniques
for the para-hyperK\" ahler manifolds (Theorem \ref{glavna1}) and
for the para-quaternionic K\" ahler manifolds (Theorem
\ref{PQKred}) and give some examples.

\section{The linear algebra  of para-quaternions}
\label{linear}
\subsection{The algebra of para-quaternions}

\noindent
In this section we will define the algebra $\PH$ of para-quaternions  and review  its basic properties.
Para-quaternions enjoy similar properties to quaternions.
The most important difference  is the existence of zero divisors in the algebra of para-quaternions.

\noindent
Both quaternions and para-quaternions are real Clifford algebras.
Let $\RR^{p,q}$ be the Euclidean space with an inner product
$\eta$ of signature  $(p,q)$,  i.e. there exists a basis $e_1, \dots , e_{p+q}$ of
$\RR^{p,q}$  such that
$\eta (e_i,e_j)=0$ for $i\ne j$, $\eta (e_i,e_i)=-1$ for $i=1, \dots , p,$ and
$\eta (e_i,e_i)=+1$ for $i=p+1, \dots ,p+q.$ The {\em real Clifford algebra
$C(p,q)$} is the universal unital real algebra which is generated by
$\RR^{p,q}$ subject to the relations
$$
v \cdot w+w \cdot v=2\eta(v,w),  \enskip v,w \in \RR^{p,q}.
$$
In particular,
$$
\HH := C(2,0) \enskip \mbox{and} \enskip \PH := C(1,1)\cong C(1,2)
$$
are algebras of quaternions and para-quaternions, respectively.
Usually, we do not write the product sign.
In other words, the algebra $\PH$ of para-quaternions is generated by unity $1$  and generators  $i,j, k$ satisfying
\begin{equation}
\label{komrel}
i^2 = -1, \enskip j^2 = 1 = k^2, \enskip ij=-ji = -k.
\end{equation}
Using  the notation $J_1 = i, J_2 = j, J_3 = k$ and constants $\e_1 := 1, \e_2 := -1 =: \e_3$, we can write the  relations (\ref{komrel}) as
$$
J_\alpha ^2 =-\e_\alpha, \enskip J_\alpha J_\beta = -\e_\gamma J_\gamma ,
$$
where $(\alpha ,\beta, \gamma)$ is a cyclic permutation of $(1, 2, 3).$

\noindent
Notice that  the relation
$$
q = a + b i + c j + d k = a + b i + j(c - d i) =: z_1(q) + jz_2(q)
$$
allows us to  identify $\PH$ with $\CC ^2$ using the map
\begin{equation}
\label{z1}
\PH \ni q \to z(q) :=
\begin{pmatrix}
z_1(q) \\
z_2(q)
\end{pmatrix}
\in \CC ^2.
\end{equation}
For a para-quaternion
$$
q = a + bi + cj + dk,
$$
with $ a, b, c, d \in \RR ,$ we define its conjugate, real and imaginary part,  by
$$
\bar q := a -  bi - cj - dk, \enskip \Re q := a,\enskip \Im q :=  bi + cj + dk,
$$
respectively. We define the scalar product on $\PH$ by
\begin{equation}
\label{scal}
\langle q, q \,' \rangle := \Re (q\bar q') = aa' + bb' - cc' - dd',
\end{equation}
where $q = a + bi + cj + dk$ and $q' = a' + b'i + c'j + d'k.$
The corresponding  square norm is multiplicative, i.e.
$$
|qq'|^2 = |q|^2|q'|^2, \enskip q, q' \in \PH.
$$
Together with the scalar product (\ref{scal}),  para-quaternions $\PH$ can be naturally identified with
$\RR ^{2,2}$; a $4$-dimensional real vector space with scalar product of  signature $(2,2).$
The commutator
$$
[q,q'] = qq' - q'q, \enskip q,q' \in \tilde H,
$$
defines a Lie algebra structure on the vector space $\PH.$  Moreover, there is a Lie algebra decomposition
$$
\PH = \RR  1 \oplus \Imag \PH ,
$$
where $\RR 1$ is the center and
$$
\Imag \PH \cong su(1,1) \cong so(2,1) \cong sl_2(\RR ) \cong sp(1, \RR )
$$
is the semisimple part.
The Lie groups corresponding to the  algebras $\PH$ and $\Imag \PH$ are
$$
\PH ^+ := \{ q \in \PH  \mid |q|^2 >0 \}, \quad \PH _1:= \{ q \in \PH  \mid |q|^2 =1 \} ,
$$
respectively.
The group $\PH_1$ of unit para-quaternions is isomorphic to $SU(1,1).$ Geometrically, $\PH _1$  is a pseudosphere $S^{2,1} \subset \RR ^{2,2},$ diffeomorphic to $S^1 \times \RR ^2.$

\subsection{Vector spaces over para-quaternions}
\label{Hn}

Consider a right module $\PH ^n \cong  \RR ^{4n}, n \geq 1$ over
the algebra $\PH$, where the multiplication of $ (h_1, \dots ,
h_n) \in \PHn$ and $q \in \PH$ is given by
$$
(h_1, \dots , h_n)q := (h_1q, \dots , h_nq).
$$
Right multiplications by $i, j, k,$  respectively, induce endomorphisms $J_1, J_2, J_3$ of $\RR ^{4n}$ satisfying
\begin{equation*}
J_\beta J_\gamma = -\e_\alpha J_\alpha, \enskip J_\alpha ^2 = -\e_\alpha \Id,
\end{equation*}
where $(\alpha, \beta , \gamma )$ is a cyclic permutation of $(1, 2, 3).$

\noindent
The map $z: \PH \to \CC ^2$ given by relation (\ref{z1}) enables us to define a real isomorphism
$z : \PHn \to \CC ^{2n}$ by the formula
\begin{equation}
\label{z2}
\PHn \ni
\begin{pmatrix}
h_1\\
\vdots \\
h_n
\end{pmatrix}
=h \to z(h) :=
\begin{pmatrix}
z(h_1)\\
\vdots \\
z(h_n)
\end{pmatrix}
\in \CC ^{2n}.
\end{equation}
If we regard $\PHn$ as a right vector space over $\CC ,$ then $z$ is a complex isomorphism.

\noindent
In $\PHn$ we define the scalar product of real signature $(2n, 2n)$ by
\begin{equation}
\label{scaln}
\langle h,  h' \rangle := \Re (h\bar h') = h_1 \bar h_1' + \dots  + h_n \bar h_n' ,
\end{equation}
where $h = (h_1, \dots , h_n),$ $h' = (h_1', \dots , h_n').$
Hence we can identify $\PHn$ with $\RR ^{2n,2n}.$
\noindent
Notice that the scalar product (\ref{scaln}) can be written in terms of the complex representation (\ref{z2}) by
\begin{equation}
\label{scalc}
\langle h,  h' \rangle = \sum_{i=1}^{n}\bigl( \Re (z_1(h_i) \overline{z_1(h_i')} - z_2(h_i) \overline{z_2(h_i')} \bigr) .
\end{equation}
It is a hermitian scalar product  of real signature $(2n, 2n)$ on $\CC ^{2n},$ represented by the following diagonal block matrix
\begin{equation}
\label{een}
E_n :=
diag(
\begin{pmatrix}
1 & 0 \\
 0 & -1
\end{pmatrix}
) \subset gl_{2n}(\CC ).
\end{equation}
With respect to the scalar product (\ref{scaln}), $J_1$ is an isometry, while $J_2$ and $J_3$ are anti-isometries of $\PHn.$
All three endomorphisms  $J_1, J_2, J_3$ are skew-symmetric with respect to the metric,
i.e.
\begin{equation*}
\langle J_\alpha \cdot , \cdot \rangle  = - \langle \cdot ,J_\alpha  \cdot \rangle  , \enskip \alpha = 1, 2, 3.
\end{equation*}

\subsection{Para-quaternionic structures in a real vector space $V$}
\label{pqinV}
\par
\begin{Definition}
\label{structures} Let $V$ be a $4n$-dimensional real vector space
with  pseudo-Rie\-mann\-ian scalar product $g$ of signature $(2n,
2n), n \geq 1.$
\begin{itemize}
\item[i)]{
A triple $(J_1, J_2, J_3)$ of endomorphisms of $V$ satisfying the relations
\begin{equation}
\label{comrel}
J_\beta J_\gamma = -\e_\alpha J_\alpha, \enskip J_\alpha ^2 = -\e_\alpha \Id,
\end{equation}
where $(\alpha, \beta , \gamma )$ is a cyclic permutation of $(1,
2, 3),$ is called a {\em para-hypercomplex structure} on $V.$
}\item[ii)]{ A subalgebra ${\cal G} \subset \End (V)$ is called a
{\em para-quaternionic structure on $V$} if there exists its basis
$J_1, J_2, J_3,$ satisfying the relations (\ref{comrel}). We say
that  the para-quaternionic structure ${\cal G}$ and the
para-hypercomplex structure  $(J_1, J_2, J_3)$ correspond to each
other. }\item[iii)]{ A para-hypercomplex structure $(J_1, J_2,
J_3)$  is called {\em hermitian with respect to $g$} if its
endomorphisms are skew-symmetric with respect to $g.$ }\item[iv)]{
A para-quaternionic structure ${\cal G}$ on $V$ is called {\em
hermitian with respect to $g$}  if some (and hence any)
corresponding  para-hypercomplex structure is hermitian with
respect to $g.$ }\end{itemize}
\end{Definition}
\noindent
Any two  bases $(J_1, J_2, J_3)$ and $(J_1' J_2', J_3')$  of ${\cal G}$ are related by  a $SO(2,1)$ transformation.
The existence of a para-hypercomplex structure $(J_1, J_2, J_3)$  allows us to equip $V$ with the structure of a right module over $\PH .$
Namely,
for $e \in V$ and $q = a + b i + c j + d k \in \PH$  we define multiplication, i.e. a right module structure on $V,$ by
$$
eq := ae + b J_1e + c J_2e + d J_3e.
$$
One can check that there always exist  vectors $e_1, \dots , e_n \in V$   such that the vectors
\begin{equation}
\label{adoptedbasis}
e_1, \dots , e_n, J_1 e_1, \dots , J_1 e_n, J_2 e_1, \dots ,J_2 e_n, J_3 e_1, \dots , J_3 e_n
\end{equation}
are linearly independent. Hence  we can make the identification
\begin{equation}
\label{indent}
V \cong \RR ^{4n} \cong \PHn .
\end{equation}
The basis (\ref{adoptedbasis}) is called the {\em basis adopted
with the para-hy\-per\-co\-mplex structure} $(J_1, J_2, J_3).$
Clearly, if the vectors  $e_1, \dots , e_n$ are orthonormal then
the adopted basis (\ref{adoptedbasis}) is a pseudo-orthonormal
basis of $V$ with respect to the scalar product $g.$ This implies
that the signature of $g$ is $(2n, 2n)$ and that this  condition
was not necessary in Definition \ref{structures}.

\noindent
To each endomorphism $J \in {\cal G}$  corresponds a nondegenerate $2$-form
$\omega _J$ on $V$ by the formula
$$
\omega _J(\cdot , \cdot ) := g( J \cdot , \cdot ) .
$$
The above correspondence is an isometry of the subspace ${\cal G} \subset \End V$ and its image in $\Lambda ^2 V.$
Let $\omega _1 := \omega _{J_1},\enskip  \omega _2 := \omega _{J_2},\enskip  \omega _3 := \omega _{J_3}.$

\noindent
One can prove that the $4$-form
\begin{equation}
\label{om} \Omega = \Omega ({\cal G}) :=\omega _1 \wedge \omega _1
- \omega _2 \wedge \omega _2 - \omega _3 \wedge \omega _3
\end{equation}
 is invariant under the action of the group $SO(2,1) $ and hence  independent of the  choice  of basis $J_1, J_2, J_3$ of ${\cal G}.$
The following lemma is a consequence of the Riemannian version by means of complexification.
\begin{Lemma}
\label{cuvaom}
The maximal subalgebra of $so(2n,2n)$ preserving the $4$-form  $\Omega ,$ given by (\ref{om}), is $sp(n, \RR) \oplus sp(1, \RR).$
\end{Lemma}

\subsection{Para-quaternionic linear maps}

\noindent
Let $V$ be a $4n$ dimensional real vector space, ${\cal G}$ a para-quaternionic structure on $V$ and ${\cal H} = (J_1, J_2, J_3)$ a corresponding para-hypercomplex structure.
\begin{Definition}
Let $L \in \End V$ be a $\RR $-linear map.

\noindent
$L$ is  {\em para-quaternionic} if for any $J \in {\cal G}$ there exists $J' \in {\cal G}$ such that $[L,J] = J'.$

\noindent
$L$ is  {\em para-hypercomplex} if $[L, J_\alpha] = 0, \enskip \alpha  = 1, 2, 3.$
\end{Definition}
\noindent
We denote the algebra  of all para-quaternionic (resp. para-hypercomplex) linear maps with respect to ${\cal G}$ (resp.  ${\cal H}$)  by  $aut({\cal G})$ (resp. $aut({\cal H})).$ By the very definitions, $aut({\cal G})$ and $aut({\cal H})$ are the normalizer and centralizer of ${\cal G}$ in $\End V.$ There is a  decomposition
$$
aut({\cal G}) = aut({\cal H}) \oplus {\cal G}.
$$
Because of the identification $V  \cong \PH  ^n$ we sometimes use
the notation $aut({\cal H}) = gl_n(\PH ).$ The algebra $gl_n(\PH
)$ can be identified with the algebra of all $n \times n $
matrices with para-quaternionic entries, acting by left
multiplication on column vectors of $\PHn .$

\noindent Let ${\cal G}$ be a hermitian para-quaternionic  and
${\cal H}$ the  corresponding para-hypercomplex structure with
respect to the  scalar product $g$ on $V.$ We define subalgebras
of $aut({\cal G})$ and $aut({\cal H})$ preserving the metric $g$
by
\begin{equation*}
sp({\cal G}) := \{ L \in aut({\cal G}) \mid g(L\cdot , \cdot ) = - g(\cdot , L\cdot ) \} ,\quad
sp({\cal H}) := aut({\cal H})  \cap sp({\cal G}).
\end{equation*}
Since ${\cal G} \subset sp({\cal G}),$ we have the decomposition
$$
sp({\cal G}) = sp({\cal H}) \oplus {\cal G}.
$$
We use the  notation $sp_n(\tilde \HH ) = sp({\cal H}).$

\begin{Lemma}
The following isomorphisms hold:\\
\centerline{ i) $gl_n(\tilde \HH ) = aut ({\cal H})  \cong
gl_{2n}(\RR ),$ \hspace{1cm} ii) $sp_n(\tilde \HH ) = sp ({\cal
H}) \cong sp(n, \RR ).$ }
\end{Lemma}
\noindent
{\bf Proof:}
First we are going to find a complex representation of the algebra $gl_n(\PH ).$
The isomorphism $z : \PHn \to \CC ^{2n}$ given by (\ref{z2}) induces the Lie algebra isomorphism $z : gl_n(\PH ) \to gl_{2n}(\CC )$ by
{\small
\begin{equation*}
\begin{pmatrix}
  a_{11} + jb_{11} & \dots &  a_{1n}+ jb_{1n}\\
  \vdots & & \vdots \\
  a_{n1} + jb_{n1} & \dots &  a_{nn}+ jb_{nn}
\end{pmatrix}
\to
\begin{pmatrix}
  \begin{array}{cc}
      a_{11} & \bar b_{11}\\
      b_{11} & \bar a_{11}
  \end{array}
  & \dots &
  \begin{array}{cc}
      a_{1n} & \bar b_{1n}\\
      b_{1n} & \bar a_{1n}
  \end{array}
  \\
  \vdots & & \vdots \\
  \begin{array}{cc}
      a_{n1} & \bar b_{n1}\\
      b_{n1} & \bar a_{n1}
  \end{array}
  & \dots &
  \begin{array}{cc}
      a_{nn} & \bar b_{nn}\\
      b_{nn} & \bar a_{nn}
  \end{array}
\end{pmatrix}
=:
\begin{bmatrix}
a_{pq} & \bar b_{pq}\\
b_{pq} & \bar a_{pq}
\end{bmatrix} .
\end{equation*}
} 

\noindent
We use the above notation for block matrices.
Hence we have the following complex representation of $gl_n(\PH ):$
$$
z(gl_n(\PH )) =
\{
\begin{bmatrix}
a_{pq} & \bar b_{pq}\\
b_{pq} & \bar a_{pq}
\end{bmatrix}
 \mid
\enskip a_{pq}, b_{pq} \in \CC,\enskip p, q = 1, .., n
\} \subset M_{2n}(\CC ).
$$
Define the map $\mu : z(gl_n(\PH ))  \to gl_{2n}(\CC )$  by
\begin{equation}
\label{mu}
\mu (Q)  = MQM^{-1},
\end{equation}
where the matrix $M\in gl_{2n}(\CC )$ is the  block diagonal matrix
$$
M = \frac{\sqrt{2}}{2}
diag(
\begin{pmatrix}
1 & i \\
i & 1
\end{pmatrix}
).
$$
One can check directly that
$$
\mu
(
\begin{bmatrix}
a_{pq} & \bar b_{pq}\\
b_{pq} & \bar a_{pq}
\end{bmatrix}
)
 =
\begin{bmatrix}
\Re a_{pq} - \Im b_{pq} & \Im a_{pq} + \Re b_{pq} \\
 \Re b_{pq} - \Im a_{pq} & \Re a_{pq} + \Im b_{pq}
\end{bmatrix}  \in gl_{2n}(\RR ).
$$
The composition  $\mu \circ z : gl_n(\PH ) \to gl_{2n}(\RR )$ is a Lie algebra isomorphism, so that statement i) is proved.

\noindent
According to the relation (\ref{scalc}), we can represent the algebra $sp_n(\PH )$ in terms of the matrix $E_n,$ given by (\ref{een}), in the following way:
\begin{align*}
(\mu \circ z)(sp_n(\PH )) &= \{ \mu (Q) \mid Q \in z(gl_n(\PH )), \mu (Q)^* \mu (E_n) = - \mu (E_n )\mu (Q) \} = \\
& = \{A \in gl_{2n}(\RR ) \mid A^TF_n = - F_n A \},
\end{align*}
where
$$
F_n =
diag(
\begin{pmatrix}
0 & 1 \\
-1 & 0
\end{pmatrix}
) = \frac{i}{2}\mu (E_n).
$$
The matrix $F_n$ represents the nondegenerate two form
$$
\omega  = e_1\wedge f_1 + \dots + e_n\wedge f_n
$$
in $\RR ^{2n},$  where $e_1, f_1, \dots , e_n, f_n$ is a basis of
$\RR ^{2n}.$ Hence $sp_n(\PH ) \cong sp(n, \RR ).$ \qed

\noindent
As an immediate consequence of the previous theorem we have the following:
\begin{Theorem}
\label{isomorphism}
Let ${\cal G}$ be a  para-quaternionic structure  and ${\cal H} = (J_1, J_2, J_3)$ a corresponding para-hypercomplex structure hermitian with respect to a scalar product $g$ on a $4n$-dimensional ($n\geq 1$) real vector space $V.$ The following isomorphisms hold:\\
\centerline{i) $aut({\cal H}) \cong gl_{2n}(\RR ),$ \hspace{1cm} $aut({\cal G}) \cong gl_{2n}(\RR )\oplus sp(1, \RR ),$}
\centerline{ii) $sp({\cal H}) \cong sp(n, \RR ),$ \hspace{1cm} $sp({\cal G}) \cong sp(n,\RR )\oplus sp(1, \RR ).$}
\end{Theorem}

\begin{Remark}
\label{case4}
The case $V \cong \RR ^4$  is covered by Theorem \ref{isomorphism}, but deserves attention. In fact, the algebra $aut({\cal G}) \cong   \RR \Id \oplus so(2,2)$ is the algebra of conformal transformations of $V.$ Moreover, by  splitting  $so(2,2) = sp(1, \RR) \oplus sp(1, \RR) $ and $sp({\cal G}) \cong {\cal G}$ we have
$$
aut ({\cal G}) = \RR \Id \oplus {\cal G} \oplus {\cal G}',
$$
where ${\cal G}'$ is another para-quaternionic structure on $V$ and $[  {\cal G}, {\cal G}'] = 0.$
\end{Remark}

\subsection{Grassman description of a para-quaternionic vector space}
\label{grassman}

\begin{Lemma}
\label{grass}
Let $V$ be a real, $4n$-dimensional vector  space.
If we represent $V$ as a tensor product
\begin{equation}
\label{tensor}
V = E \otimes H,
\end{equation}
where $E \cong \RR ^{2n}$, $ H \cong \RR ^2,$ then the endomorphisms
\begin{equation}
\label{J123}
J_1 = \Id \otimes
\begin{pmatrix}
0 & -1\\
1 & 0
\end{pmatrix}
, \quad
J_2 = \Id \otimes
\begin{pmatrix}
0 & 1\\
1 & 0
\end{pmatrix}
, \quad
J_3 = \Id \otimes
\begin{pmatrix}
-1 & 0\\
0 & 1
\end{pmatrix}
\end{equation}
of $V$ form  a para-hyperK\" ahler structure on $V.$ Conversely,
for any para-hyperK\" ahler structure $(J_1, J_2, J_2)$ on $V$
there exists a decomposition (\ref{tensor}) of $V$ and a basis of
$H$ in which endomorphisms $J_1, J_2, J_3$ are of the form
(\ref{J123}).
\end{Lemma}
\noindent {\bf Proof:} The first statement is obvious. To check
the second let ${\cal H} = (J_1, J_2, J_3)$ be a para-hypercomplex
structure on $V$ and let $H :=\RR \langle \Id - J_3, J_1 + J_2
\rangle .$ The space $H$ is a subalgebra of $\RR \langle \Id, J_1,
J_2, J_3 \rangle \cong \PH .$ Let
$$
E := \RR \langle e_1, \dots e_n, J_2 e_1, \dots , J_2 e_n \rangle ,
$$
for  vectors  $e_1, \dots , e_n \in V,$ such that (\ref{adoptedbasis}) is a basis of $V.$
Then the map $\phi : E \otimes H  \to V$ defined by
$$
\phi (e\otimes h ) := h(e), \enskip e\otimes h \in E\otimes H = V ,
$$
is a vector space isomorphism.
It is easy to check that the endomorphisms $J_1, J_2, J_3$ have the form (\ref{J123}) in the basis $h_1 =
 \Id - J_3,$ $h_2 = J_1 + J_2$ of $H.$\qed

\noindent It is clear that  Lemma \ref{grass} holds for the
para-quat\-er\-ni\-onic stru\-ctu\-re ${\cal G} =  \RR  \langle
J_1, J_2, J_3\rangle .$

\noindent
The algebra $aut ({\cal G}) \cong  gl_{2n}{\RR} \oplus sp(1, \RR )$ acts irreducibly on $V  = E \otimes H$ by
\begin{equation}
\label{action}
(A, J) (e\otimes h) = Ae \otimes h + e \otimes Jh,
\end{equation}
for any $(A, J) \in gl_{2n}{\RR} \oplus sl(1, \RR ),$ $e\otimes h \in V.$
The algebra $sp({\cal G}) = sp(n, \RR ) \oplus sp(1, \RR )$
acts on  $V= E\otimes H$  by the formula (\ref{action}) preserving some nondegenerate two forms $\omega ^E$ and $\omega ^H$ on $E$ and $H,$ respectively.
The scalar product
\begin{equation}
\label{geom}
g = \omega ^E \otimes \omega ^H
\end{equation}
 on $V$ of  signature $(2n, 2n)$ is invariant under the action of $sp({\cal G}).$
The para-quaternionic structure ${\cal G} =  \RR  \langle J_1, J_2, J_3\rangle $ is hermitian with respect to the scalar product (\ref{geom}).

\noindent
Conversely, for any para-quaternionic structure ${\cal G},$ hermitian with respect to a given metric $g$ on $V,$ there exist forms $\omega ^E$ and $\omega ^H$ such that relation (\ref{geom}) holds.

\noindent
We conclude this section  with the following interesting  lemma:
\begin{Lemma}
\label{dist1}
Let ${\cal G}$  and ${\cal H} = (J_1, J_2, J_3)$ be a para-quaternionic and corresponding para-hypercomplex structure  hermitian with respect to a scalar product $g$ on a vector space $V.$ Then
\begin{itemize}
\item[i)]{
The algebra $sp({\cal G}) = sp(n, \RR ) \oplus sp(1, \RR )$ acts on $V$ irreducibly.
}\item[ii)]{
The algebra $sp({\cal H}) = sp(n, \RR ) $ acts on $V$ reducibly. Moreover, there exist a $S^1$ family of isotropic $2n$-dimensional subspaces  $V_\phi,$ $\phi \in S^1,$ invariant with respect to the action of $sp({\cal H}).$
}\end{itemize}
\end{Lemma}
\noindent
{\bf Proof:} The first statement   is obvious since  the action (\ref{action}) is a tensor product of irreducible actions.
To prove ii), let $h_1, h_2$ be a basis of $H.$  For any $\phi \in [0, \pi ]$ the $2n$-dimensional subspace
$$
V_\phi := \{e\otimes ( h_1\cos \phi + h_2\sin \phi ) \mid e \in E \}
$$
is $sp({\cal H})$ invariant and isotropic. \qed

\section{Para-quaternionic structures on manifolds}
\label{manifolds}
\subsection{Para-hyperK\" ahler manifolds}

\begin{Definition}
\begin{itemize}
\item[i)]{
A  pseudo-Riemannian manifold $M =(M^{4n}, g), n \geq 1$ is called
{\em almost hermitian para-hy\-per\-co\-mplex} if there exist
three global sections $J_1,J_2, J_3$ of $\End(TM)$ with the
following property: for  each point $p \in M$ the  triple $(J_1,
J_2, J_3)$ is a hermitian para-hypercomplex structure with respect
to the metric $g$ of  $T_pM$. }\item[ii)]{ \noindent An almost
hermitian para-hypercomplex manifold $M =(M^{4n}, g), n \geq 1$ is
{\em para-hyperK\" ahler} if $\nabla J_\alpha = 0, \enskip \alpha
= 1, 2, 3,$ and $\nabla$ is the Levi-Civita connection with
respect to the metric $g.$ }\end{itemize}
\end{Definition}
\noindent
Note that the definition implies that the signature of the metric $g$ is $(2n, 2n).$

\begin{Theorem}
The following conditions are equivalent for an almost hermitian para-hypercomplex manifold $M =(M^{4n}, g, (J_1, J_2, J_3))$ ($\nabla$ denotes the Levi-Civita connection):
\begin{itemize}
\item[i)]{
$\nabla  J_\alpha  = 0, \enskip \alpha = 1,2, 3,$
}\item[ii)]{
$\nabla  \omega _\alpha  = 0, \enskip \alpha = 1,2, 3,$
}\item[iii)]{
$d \omega _\alpha  = 0, \enskip \alpha = 1,2, 3.$
}\end{itemize}
\end{Theorem}
The proof is similar as in the Riemannian case.\qed

\begin{Remark}
The para-hyperK\" ahler manifold is K\" ahler since $J_1$ is an integrable complex structure and $d\omega _1 = 0.$ The structures $ J_2$ and $J_3$ are not complex, but they are product structures.
\end{Remark}
\noindent
The following characterization of para-hyperK\" ahler manifolds follows from the definition of  the algebra $sp({\cal H}) \cong sp(n, \RR)$ and the discussion in  Section \ref{grassman}.
\begin{Theorem}
\label{charaPHK}
A pseudo-Riemannian manifold $(M^{4n}, g)$ is para-hyperK\" ahler  if and only if $hol(M ) \subset sp(n, \RR) .$
\end{Theorem}

\noindent According to the results of Section \ref{grassman} we
can decompose the tangent space of a para-hyperK\" ahler manifold
$M$ as $TM = E \otimes H,$ where $E$ and $H$ are parallel
$2n$-dimensional and $2$-dimensional  vector bundles,
respectively. The metric $g$ of the manifold $M$ can be written in
the form
$$
g = \omega ^E \otimes \omega ^H
$$
where $\omega ^E$ and $\omega ^H$ are certain nondegenerate forms on $E$ and $H,$ respectively.
The Levi-Civita connection $\nabla$ is $sp(n, \RR)\otimes  T^*M$ valued and is given by
$$
\nabla (e\otimes h ) = \nabla ^Ee \otimes h + e \otimes \nabla^H h
$$
for any $e\otimes h \in \Gamma (TM).$ The conditions that $ \nabla$ is Levi Civita and commutes with  the hypercomplex structure are equivalent to
$$
\nabla ^H \omega ^H =0 = \nabla ^E \omega ^E, \quad \nabla ^H \mbox{ is flat.}
$$
An interesting consequence of the flatness of $ \nabla ^H$  and Lemma \ref{dist1} b) is the following:
\begin{Theorem}
Let $M = (M^{4n}, g)$  be a para-hyperK\" ahler manifold.
There exists a $S^1$ family $V_\phi, \enskip \phi \in S^1,$ of isotropic, parallel, $2n$-dimensional distributions of $TM.$
\end{Theorem}
\noindent
However, since the distributions $V_\phi$ are isotropic, this does not imply decomposability of the manifold $M.$

\begin{Example}
\label{phn}
The para-quaternionic vector space $\PHn$ with the metric $g$ given by formula (\ref{scaln}) is a para-hyperK\" ahler manifold.
The required endomorphisms  $J_1, J_2, J_3$ of $T\PHn \cong \PHn$ can be  defined as right multiplications by $i, j$ and  $k,$ respectively.
\end{Example}
\noindent Using the method of para-hyperK\" ahler reduction
 from this flat example Section \ref{hkred} we will
obtain another example of para-hyperK\" ahler manifold.

\begin{Example}
All Riemannian homogenous hyperK\" ahler manifolds are flat. In
particular, all Riemannian symmetric hyperK\" ahler spaces are
flat.  The indefinite hyperK\" ahler symmetric spaces are not
necessarily flat. They are classified in \cite{AC}. Indefinite
hyperK\" ahler and para-hyperK\" ahler symmetric spaces have the
same complexification. Accordingly, classification of complex
hyperK\" ahler symmetric spaces from \cite{AC} lead to
classification of para-hyperK\" ahler symmetric spaces in
\cite{ABCV}. It is interesting fact that para-hyperK\" ahler
symmetric spaces exist in all dimensions $4n$ unlike the
indefinite hyperK\" ahler symmetric spaces which exist only in
dimensions $8n.$ For example, there are up to the isomorphism only
two (mutually anti-isometric) non-flat para-hyperK\" ahler
symmetric spaces given as follows.

 \noindent
 Let us take the matrix $A$
$$
A =\left[ \begin{matrix}
 0 & -1/2 & 1/2 & 0 \\
1/2 &0 & 0 & 1/2 \\
1/2 &0& 0 & 1/2 \\
 0 & 1/2& -1/2 & 0
\end{matrix}
\right]
$$
Define a subalgebra $f= \RR A \subset so(2,2)$ acting on the
vector space $m= \RR ^{2,2}$ from the left. Let $E_1, E_2, E_3,
E_4$ be a pseudo-orthonormal basis of $m.$ Define the  Lie
algebras $g_\pm =  m + f$ by the nonzero commutators
\begin{align*}
[A, M] &= A(M), \enskip M \in m,\\
[E_1, E_2]=[E_3, E_1]&=[E_4, E_2]=[E_3, E_4]= \pm A
\end{align*}
The symmetric spaces corresponding to the symmetric decompositions
$g_\pm = m + f$ are required para-hyperK\" ahler spaces. These
spaces appeared in \cite{BBR} as examples of rank $2$ Osserman
symmetric spaces. Note that the algebras $g_\pm$ are solvable.
\end{Example}

\subsection{Para-quaternionic K\" ahler manifolds}

\begin{Definition}
A pseudo-Riemannian manifold $M =(M^{4n}, g), \enskip n \geq 2,$
is {\em almost hermitian  para-quaternionic} if there exists a
subbundle  ${\cal G}$ of $\End(TM)$ with the  following property:
for each point $p \in M$ the fiber ${\cal G}_p$ of ${\cal G}$ is a
hermitian para-quaternionic structure of $T_pM$ with respect to
the metric $g$.
\end{Definition}
\noindent
We can locally choose a para-hypercomplex  structure  $(J_1, J_2, J_3)$ which is a basis of ${\cal G}.$ Hence any almost hermitian para-quaternionic manifold is locally almost hermitian para-hypercomplex.

\noindent
As in Section \ref{pqinV} we can locally define nondegenerate two forms $\omega _1, \omega _2, \omega _3.$
However the $4$-form
$$
\Omega = \Omega ({\cal G})= \omega _1 \wedge \omega _1 - \omega _2 \wedge \omega _2 - \omega _3 \wedge \omega _3
$$
is defined globally on $M.$

\begin{Definition}
\begin{itemize}
\item[i)]{
An almost hermitian para-quaternionic manifold $M=(M^{4n}, g),$ $n \geq 2,$ is called {\em para-quaternionic K\" ahler}  if $\nabla \Omega = 0$ and $M$ is not para-hyperK\" ahler.
}\item[ii)]{
A pseudo-Riemannian $4$-dimensional manifold $(M^4, g)$ with a metric $g$ of  signature $(2,2)$ is called {\em para-quaternionic K\" ahler} if it is self-dual and Einstein.
}\end{itemize}
\end{Definition}
\noindent
Clearly, the condition $\nabla \Omega = 0$ implies $d\Omega = 0$ which has strong consequences on the topology of $M.$

\noindent According to Lemma \ref{cuvaom}  and Theorem
\ref{charaPHK} for a para-quaternionic K\" ahler manifold  we have
$hol(M) \subset sp(n, \RR )\oplus sp(1, \RR),$ $hol(M) \not
\subset sp(n, \RR ).$ The converse statement follows from the
results of Section \ref{grassman}, so we have the following
characterization.

\begin{Theorem}
\label{charaPQK}
A pseudo-Riemannian manifold $M = (M^{4n}, g), \enskip n \geq 2,$ is para-qua\-te\-rni\-onic K\" ahler if and only if
$$
hol(M) \subset sp(n, \RR )\oplus sp(1, \RR),\enskip  hol(M) \not \subset sp(n, \RR ).
$$
\end{Theorem}
\noindent Note that for a $4$-dimensional manifold $M$ with a
metric $g$ of signature $(2,2)$ the condition  $hol(M) \subset
so(2, 2) = sp(1, \RR ) \oplus sp (1, \RR )$ is always satisfied.
On the other hand, the condition $\nabla \Omega = 0$ is satisfied
since $\Omega$ is a $4$-form. These were the reasons for the
separate definition of a $4$-dimensional para-quaternionic K\"
ahler manifold.

\noindent
Similarly to the para-hyperK\" ahler case we can write  $TM = E \otimes H$ where $E$ and $H$ are certain $2n$-dimensional and $2$-dimensional  vector bundles, respectively.
If we write the metric in the form
$$
g = \omega ^E \otimes \omega ^H
$$
for some $2$-forms $\omega ^E$ and $\omega ^H$ then the
Levi-Civita connection $\nabla$ is given by
$$
\nabla (e\otimes h ) = \nabla ^Ee \otimes h + e \otimes \nabla^H h
$$
for any $e\otimes h \in \Gamma (TM).$ The condition that $ \nabla$ is Levi-Civita is equivalent to
$$
\nabla ^H \omega ^H =0 = \nabla ^E \omega ^E.
$$
The connection $\nabla$ preserves the subbundle ${\cal G} \subset \End (TM).$ One can easily check that the structure equations are
$$
\nabla J_\alpha = \e _\alpha (\mu _\beta J_\gamma - \mu _\gamma J_\beta ),
$$
where $(\alpha, \beta, \gamma )$ is a cyclic permutation of $(1,2,3)$ and $\mu _1, \mu _2, \mu _3$ are structure $1$-forms.

\begin{Example}
We shall show that the quotient $sl_{n+2}(\RR )/ sl_{n}(\RR )
\oplus sp(1, \RR)$ is a para-quaternionic Ka\" hler,
pseudo-Riemannian symmetric space. Its holonomy  $sl_{n}(\RR )
\oplus sp(1, \RR) $ $\subset sp(n, \RR) \oplus sp(1, \RR)$ is not
full. An example of a para-quaternionic symmetric space with full
holonomy is given in Section \ref{projspace}.

\noindent
Consider the symmetric decomposition   $sl_{n+2}(\RR ) = m + f$ of the algebra $sl_{n+2}(\RR ),$ where
$$
f = \{
\begin{pmatrix}
A & 0\\
0 & B
\end{pmatrix} \mid A \in sl_2(\RR ), B \in sl_n(\RR )
\} \cong  sl_n(\RR ) \oplus sp(1,\RR ) .
$$
One can easily check that the adjoint action of the algebra $sp(1,
\RR )$ defines a para-quaternionic structure ${\cal G}$ on the
vector space $m \cong \RR^{4n}$ invariant under holonomy.  Hence
the pseudo-Riemannian symmetric space $sl_{n+2}(\RR )/ (sl_{n}(\RR
)\oplus sp(1, \RR))$ is para-quaternionic K\" ahler. One can  find
a basis of ${\cal G}$ and then decompose  $m = E\otimes H$ using
the proof of Theorem \ref{grass}.
\end{Example}

\subsection{The para-quaternionic projective space $\PH P^n$}
\label{projspace}

Consider the vector space $\PH ^{n+1}$ with a natural scalar product $\langle \cdot, \cdot \rangle $ of signature $(2(n+1), 2(n+1))$ defined in Section \ref{Hn}.
Denote by
$$
S =  S^{2(n+1), 2n+1} \subset \PH ^{n+1}
$$
the pseudosphere of unit
vectors in $\PH ^{n+1}$  with the induced metric. Its signature is $(2(n+1), 2n+1).$

\noindent
The group $\H1 $ of unit para-quaternionic numbers acts freely and isometrically on
$S$ from the right by the formula
$$
xq = (x_1,\dots ,x_{n+1})q := (x_1q,\dots ,x_{n+1}q),
$$
$x=(x_1,\dots ,x_{n+1}) \in S, \enskip q \in \H1.$
The para-quaternionic projective space $\PHPn$ is defined as the orbit space
$$
\PHPn := S/\H1 .
$$
Denote by
$$
\pi : S \to \PHPn, \quad \pi ((x_1,\dots ,x_{n+1}) ):=[x_1,\dots ,x_{n+1}]
$$
the natural projection.

\noindent Since the orbits of $\H1$ are nondegenerate of signature
$(2,1)$ the metric $\langle \cdot, \cdot \rangle$ pushes down to
the pseudo-Riemannian metric $g$ of signature $(2n,2n)$ on $\PHPn
.$ The manifold $\PHPn$ has the  constant para-quaternionic
sectional curvature (see \cite{B}).

\noindent
Let us show that $\PHPn$ is a para-quaternionic K\" ahler manifold.
The vertical subspace $T_x^vS$  of $T_xS, \enskip  x \in S,$ with respect to the action of  $\H1$ is
$$
T_x^vS :=Ker(d_x\pi ) = \RR \langle xi, xj, xk \rangle   \subset \PH ^{n+1}.
$$
\noindent
Moreover, there is an orthogonal  decomposition
$$
T_xS = T_x^vS \oplus T_x^hS
$$
of $T_xS,$ $x \in S,$  invariant to the action of $\H1 .$
The right multiplication by $h \in \Imag \PH$ preserves $T_x^hS.$
Let
$$
s : \PHPn \to S
$$
be a horizontal section.
For any $h \in \Imag \PH$ define the  endomorphism $J_h$ of $T\PHPn$ by
$$
J_h(X):=\pi_*(s_*(X)h),
$$
\noindent
where $\enskip X \in T\PHPn.$
\noindent
Since the endomorphisms $J_h$ and $J_h'$ corresponding to different sections $s$ and $s' = sq, \enskip q \in \H1 ,$ are $Ad(q)$ related, they generate a  global subbundle ${\cal G}$ of $\End (T\PHPn ).$
One can easily check that ${\cal G}$ is an almost  para-quaternionic K\" ahler structure on $\PHPn $ and that the corresponding $4$-form $\Omega ({\cal G})$ is parallel.
Hence $\PHPn$ is a para-quaternionic K\" ahler manifold.

\begin{Theorem}
\label{hom}
The para-quaternionic projective space $\PHPn$ is the homogenous space
\begin{equation}
\label{homogenous}
\PHPn = Sp(n+1, \RR ) / Sp(1,\RR )\cdot  Sp(n,\RR )
\end{equation}
where $Sp(1,\RR )\cdot  Sp(n,\RR ) = (Sp(1, \RR )\times   Sp(n,\RR ))/Z_2$ and $Sp(n, \RR)$ denotes the (connected) group of symplectic transformations.
\end{Theorem}
\noindent
{\bf Proof:}
The map $\mu$ defined by formula (\ref{mu}) is a Lie group isomorphism and
$$
Sp(n+1, \RR ) \cong \mu ^{-1} (Sp(n+1, \RR )) = \{
\begin{pmatrix}
A & B \\
B & \bar A
\end{pmatrix}
\in SU(n+1,n+1) \}
$$
holds. Using the representation $z : \PH ^{n+1} \to \CC ^{2(n+1)}$ given by (\ref{z2}), we obtain an isometric action of  $\mu ^{-1} (Sp(n+1, \RR ))$ on the pseudosphere $S \subset  \PH ^{n+1}.$ This action commutes with the right multiplication by the para-quaternionic numbers defining the space $\PHPn$ and hence the group $\mu ^{-1} (Sp(n+1, \RR ))$ acts on $\PHPn.$
The induced action is isometric and preserves the para-quaternionic-structure ${\cal G}$ on $\PHPn.$
To show its transitivity we will show that the point
$$
o = (1, 0, \dots , 0) + j(0,0, \dots , 0) \in S \subset \PH ^{n+1}
$$
can be mapped to any point $ a_1 + jb_1 \in S$ by a matrix  $M \in \mu ^{-1} (Sp(n+1, \RR ))$.

\noindent
Notice that $|a_1 + jb_1|^2 = 1$ imply
$$
|\bar b_1 + j\bar a_1|^2 = -1,\enskip  \langle a_1 + jb_1, \bar b_1 + j\bar a_1 \rangle = 0
$$
and $\langle a_1 + jb_1, \bar a_2 + j\bar b_2 \rangle = 0$, $|a_2 + jb_2|^2 = 1$ implies
$$
\langle \bar b_1 + j\bar a_1 , \bar b_2 + j\bar a_2 \rangle = 0.
$$
Hence we can construct a  matrix $M \in \mu^{-1}(Sp(n+1,\RR )) \subset SU(n+1, n+1)$ having the vectors
$$
(z(a_1 + jb_1) , z(\bar b_1 + j\bar a_1), \dots , z(a_{n+1} + jb_{n+1}), z(\bar b_{n+1} + j\bar a_{n+1}))
$$
as columns and mapping the point $o \in S$ to a point $a_1 + jb_1 \in S.$
One can easily check that the stabilizer of the point $[o] = \pi (o) \in \PHPn$ is isomorphic to
$Sp(1,\RR )\cdot  Sp(n,\RR ).$
\qed

\begin{Theorem}
The space $\PHPn$ is a para-quaternionic K\" ahler, pseudo-Riemannian symmetric space
$$
\PHPn = sp(n+1, \RR)/(sp(n, \RR) \oplus sp(1, \RR)).
$$
 Its curvature tensor is given by the formula
\begin{equation}
\label{curvature}
\begin{split}
R(A, B)C  & =  g(B, C)A - g(A,C)B +
\sum _\alpha \e_\alpha \bigl( g(J_\alpha B, C) J_\alpha A -  \\
 & - g(J_\alpha A, C) J_\alpha B \bigr)
 -  2\sum _\alpha \e_\alpha g(J_\alpha A, B)J_\alpha C.
\end{split}
\end{equation}
\end{Theorem}
\noindent
{\bf Proof:}
Since the calculations are simple unlike the notation, we give only the idea.
The notation is the same as in the  proof of Theorem \ref{hom}.
\noindent
First one  should find the symmetric algebra decomposition
\begin{equation}
\label{split1}
sp(n+1, \RR ) \cong \mu^{-1} (sp(n+1, \RR)) = m' + f'
\end{equation}
by interpreting the proof of the Theorem \ref{hom}.
Then one should find expressions for the  metric $g'$ and endomorphisms $J_1', J_2', J_3'$ in terms of matrices from $m'$ (representing tangent vectors).
\noindent
The next step is to carry  the decomposition (\ref{split1}) together with all corresponding structures into the decomposition
\begin{equation}
\label{split2} sp(n+1, \RR) = m + f, \quad f = sp(n, \RR ) \oplus
sp(1, \RR ),
\end{equation}
by the isomorphism $\mu .$
 Using the expression $R(A,B)C = [[A,B], C]$ for the curvature of a symmetric space one easily checks the relation (\ref{curvature}).\qed

\begin{Remark}
Let $e^1, f^1, \dots , e^{n+1}, f^{n+1}$ be a  basis and $\omega=
e^1\wedge f^1 +  \dots + e^{n+1}\wedge f^{n+1}$ be a nondegenerate
$2$-form on $\RR^{2(n+1)}.$ One can identify $sp(n+1, \RR)$ and $
S^2\RR ^{2(n+1)} $ by
$$
(a \vee b)(x)= \omega (a,x)b + \omega (b,x)a, \enskip x \in
\RR^{2(n+1)},
$$
where $a \vee b \in S^2\RR ^{2(n+1)}$ denotes a symmetric product of vectors.
The decomposition (\ref{split2}) is obtained after stabilizing $\RR \langle e_0, f_0 \rangle$ under that action.
\end{Remark}

\noindent
Concerning the topology of para-quaternionic projective spaces the following is known (see \cite{BV1}).
\begin{Theorem}
The space $\PHPn$ is homotopically equivalent to the complex projective space $\CC P^n.$
\end{Theorem}

\section{Space  of curvature tensors of the para-quaternionic type}
\label{decomposition} \noindent In this we shall just state the
basic results. For details the reader is referred to \cite{AM}
where a similar decomposition has been done in the quaternionic
case.

\begin{Definition}
Let $V$ be a real vector space  and let ${\cal L}$ be  any subalgebra of $\End  (V).$
The {\em space ${\cal R}({\cal L})$ of curvature tensors of type ${\cal L}$} is the space of ${\cal L}$ valued $2$-forms $R$ on $V$ satisfying the Bianchi identity
$$
R(X, Y)Z + R(Y, Z)X + R(Z, X)Y = 0,
$$
for all $X,Y,Z \in V.$
\end{Definition}
\noindent
The space ${\cal R}({\cal L})$  is an ${\cal L}$ module.

\noindent Let ${\cal G} \subset \End (V)$  be a para-quaternionic
structure on a $4n$-dimensional real vector space $V$ and ${\cal
H} = (J_1, J_2, J_3)$  a corresponding para-hypercomplex
structure. To describe the decomposition of the module ${\cal
R}({\cal L})$ on irreducible submodules in case ${\cal L} \in \{
aut({\cal G}), aut ({\cal H}) \}$ we need some additional
notation.

\noindent
We say that a bilinear form $B \in Bil (V)$  is {\em hermitian with respect to a para-quaternionic structure} ${\cal G}$ if
$$
B(JX, Y) = - B(X, JY),
$$
for all $J \in {\cal G}$, $X, Y \in V.$ Denote by $Bil_{\cal
G}(V)$ the set of all bilinear forms on the space $V$ with respect
to ${\cal G}$ and by $\Pi_{\cal G} : Bil (V) \to Bil _{\cal G}(V)$
the projector
$$
\Pi_{\cal G}(B)(X, Y) := \frac{1}{4} \bigl( B(X, Y) + \sum_{\alpha =1}^3\e_\alpha B(J_\alpha X, J_\alpha Y) \bigr) ,
$$
$B \in Bil(V), X, Y \in V,$ which is well defined, i.e. independent of the basis $J_1, J_2, J_3$ of ${\cal G}.$
The set of bilinear forms decomposes as follows:
$$
Bil (V) = S^2_{\cal G} + \Lambda ^2 _{\cal G} + S^2_{mix} + \Lambda ^2 _{mix} ,
$$
where $S^2_{\cal G}$ and $\Lambda ^2 _{\cal G}$ are the  sets of symmetric and antisymmetric forms in $Bil_{\cal G} (V),$ respectively, and  $S^2_{mix} := (\Id - \Pi _{\cal G})(S^2_{\cal G})$ and $\Lambda ^2_{mix} := (\Id - \Pi _{\cal G})(\Lambda ^2 _{\cal G})$ their complements.

\noindent
For any $B \in Bil (V)$ define
\begin{align}
\label{curvB} R^B(X, Y)Z   & :=(B(Y,X) - B(X, Y))Z  + \\ \nonumber
 + \sum_\alpha \e _\alpha \bigl( B(X, J_\alpha Y) & - B(Y, J_\alpha X))Z +
 B(X, J_\alpha Z) J_\alpha Y - B(Y, J_\alpha Z) J_\alpha X \bigr) ,
\end{align}
for all $X, Y, Z \in V.$
One can check that the map
$$
\phi : Bil (V) \to {\cal R}(aut ({\cal G})),  \quad \phi (B) := R^B
$$
is a well defined monomorphism.
Its image
$$
{\cal R}^{Bil} : = \phi(Bil(V)) \subset {\cal R}(aut ({\cal G}))
$$
is an $aut ({\cal G})$ module.
\begin{Remark}
Let the para-quaternionic structure ${\cal G}$ be hermitian with respect to the  scalar product $g$ on $V.$
Comparing the expression (\ref{curvB}) for $B=g$ and the expression (\ref{curvature}) for the curvature at any point of the para-quaternionic projective space $\PHPn,$  we can see that  $R^g = R_{\PHPn }.$
\end{Remark}
\noindent
For a curvature tensor $R$ we define its Ricci tensor $Ric(R)$ by
$$
Ric(R)(Y, Z) := \Tr (X\mapsto R(X, Y)Z),
$$
$X, Y, Z \in V.$ One can check that the Ricci map $Ric : {\cal
R}^{Bil} \to Bil (V)$ is a monomorphism. Hence a curvature tensor
$R \in {\cal R}(aut({\cal G}))$ has an unique decomposition $R =
W(R) + Ric(R)$ with respect to the decomposition
$$
{\cal R}(aut({\cal G})) = {\cal W} \oplus  {\cal R}^{Bil}, \quad {\cal W} := Ker(Ric).
$$
 $W(R)$ is called {\em Weil part} and $Ric(R)$ is called {\em Ricci part} of the tensor $R.$

\noindent
Now we can state the main theorem:
\begin{Theorem}
\label{decompGL}
Let ${\cal G}$  be a para-quaternionic structure on a $4n$-dimensional vector space $V$ and ${\cal H} = (J_1, J_2, J_3)$  a corresponding para-hypercomplex structure.
\begin{itemize}
\item[i)]{
For $n \geq 2$ the module  ${\cal R}(aut ({\cal G}))$ is the sum of the following irreducible submodules
\begin{equation*}
{\cal R}(aut({\cal G})) = {\cal W} + {\cal R}(S^2_h) + {\cal R}(S^2_{mix}) + {\cal R}(\Lambda ^2_h) + {\cal R}(\Lambda ^2_{mix}),
\end{equation*}
where ${\cal W}$ is the kernel of the Ricci map and ${\cal R}(S^2_h),$ ${\cal R}(S^2_{mix}),$ ${\cal R}(\Lambda ^2_h),$ ${\cal R}(\Lambda ^2_{mix})$ are images  of $S^2_{\cal G},$   $\Lambda ^2 _{\cal G},$ $S^2_{mix}, $ and $\Lambda ^2 _{mix},$ respectively, by the map $\phi .$

\noindent
  For $n =1$ the above decomposition   holds but the module ${\cal W}$ is decomposable as
$$
{\cal W} = {\cal W}_+ + {\cal W}_- = {\cal R}({\cal G})+ {\cal R}({\cal G}'),
$$
where ${\cal G}'$ is  para-quaternionic structure on $V$ described
in Remark \ref{case4}. }\item[ii)]{ ${\cal R}(aut({\cal H}))  =
{\cal W}  + {\cal R}(\Lambda ^2_h), \enskip n \geq 2,$ \hspace{1
em} ${\cal R}(aut({\cal H}))  = {\cal W}_-  + {\cal R}(\Lambda
^2_h), \enskip n =1 .$ }\end{itemize}
\end{Theorem}
\begin{Theorem}
\label{decompSP} Let ${\cal G}$  be a para-quaternionic hermitian
structure  with respect to the scalar product $g$ on a
$4n$-dimensional,  vector space $V,$ $n \geq 2$ and let ${\cal H}
= (J_1, J_2, J_3)$  be a corresponding para-hypercomplex
structure.
\begin{itemize}
\item[i)]{ Any curvature tensor $R \in {\cal R}(sp({\cal G}))$  is Einstein i.e.
$$
Ric(R) = \frac{K(R)}{4n}g,
$$
holds, where $K(R)$ is the scalar curvature of $R.$ }\item[ii)]{
${\cal R}(sp({\cal G})) = \RR R_{\PHPn } + {\cal R}(sp({\cal H}))$
where ${\cal R}(sp({\cal H})) \subset {\cal W}$ is an irreducible
submodule and $R_{\PHPn }$ is the curvature of the
para-quaternionic projective space. }\end{itemize}
\end{Theorem}

\noindent
We have the following important consequences  concerning para-quaternionic K\" ahler manifolds.
\begin{Theorem}

\begin{itemize}
\item[i)]{A para-quaternionic K\" ahler manifold is an Einstein manifold.
}\item[ii)]{ A para-quaternionic K\" ahler manifold of dimension $4n, \enskip n \geq 2,$ with zero scalar curvature  is a locally para-hyperK\" ahler manifold.
}\end{itemize}
\end{Theorem}

\noindent
The following lemma is  important both for the  proofs of the previous theorems and for various applications.
\begin{Lemma}
Let ${\cal G}$  be a para-quaternionic hermitian structure with respect to the  scalar product $g$ on a $4n$-dimensional,  vector space $V$ and ${\cal H} = (J_1, J_2, J_3)$  a corresponding para-hypercomplex structure.
\begin{itemize}
\item[i)]{A curvature tensor $R \in {\cal R}(\End (V))$ belongs to ${\cal R}(aut({\cal G}))$ if and only if
$$
[R, J_\alpha ] = \frac{\e_\alpha}{2n} \bigl( \Tr (J_\gamma R)J_\beta -\Tr (J_\beta R)J_\gamma  \bigr),
$$
where $(\alpha, \beta, \gamma )$ is a cyclic permutation of $(1, 2, 3).$
}\item[ii)]{
Tensors belonging to the module ${\cal W}$  are traceless.
}\item[iii)]{
The following commutator relations hold:
$$
[{\cal W}, {\cal G}] = 0, \enskip n \geq 2 \quad \mbox{ and } \quad [{\cal W}, {\cal G}'] = 0, \enskip n = 1,
$$
where the para-quaternionic structure ${\cal G}'$ is described in Remark \ref{case4}.
}
\end{itemize}
\end{Lemma}

\section{Reduction techniques}
\label{reduction}
\subsection{Para-hyperK\" ahler reduction}
\label{hkred}

Let $M = ( M^{4n}, g, (J_1, J_2, J_3))$ be a para-hyperK\" ahler manifold and let $\omega _1, \omega _2, \omega _3$ be the global nondegenerate $2$-forms associated with $J_1, J_2, J_3.$
Let a  Lie  group  $G$ acts  freely on $M$ by isometries and preserves para-hyperK\" ahler structure i.e.
$g^*J_\alpha  = J_\alpha , \enskip \alpha= 1, 2, 3,$ holds for a $g \in G.$
Let $V$ be a Killing vector field on $M$ generated by an element $V^* \in {\cal g} = Lie (G).$

\noindent
If it exists, the map $f = (f_1, f_2, f_3) : M \to {\cal g}^*\otimes {\cal g}^*\otimes {\cal g}^*$
such that
$$
\noindent
d(f_\alpha (V)) = \omega _\alpha (V^*,\cdot ) ,
$$
$\enskip \alpha = 1, 2, 3,$
holds for any $V \in {\cal g}$ at any point of the manifold $M,$ is called the {\em moment map for the action $G$ on $M$}.

\begin{Theorem}[Para-hyperK\" ahler reduction]
\label{glavna1} Let a Lie group $G$  act  freely and isometrically
on a para-hyperK\" ahler manifold $M = ( M^{4n}, g, (J_1, J_2,
J_3))$  and preserve its para-hyperK\" ahler structure. Let $f =
(f_1, f_2, f_3) : M \to {\cal g}^*\otimes {\cal g}^*\otimes {\cal
g}^*$ be the corresponding equivariant moment map. Suppose that
$\xi \in {\cal g}^*\otimes {\cal g}^*\otimes {\cal g}^*$ is such
that
$$
{\cal K}_\xi : = f^{-1}(\xi )
$$
(on which $G$ acts by isometries) is a smooth submanifold of $M.$
Suppose that the quotient
$$
M_\xi :=  {\cal K}_\xi /G
$$
 has a smooth manifold structure for which the projection $\pi _\xi : {\cal K} _\xi \to M_\xi$ is a smooth pseudo-Riemannian submersion and $g_\xi$ is the induced metric on $M_\xi .$ Then  the pseudo-Riemannian manifold $(M_\xi , g_\xi )$ is a para-hyperK\" ahler manifold with respect to the para-hyperK\" ahler structure  obtained by projection $\pi$ from the structure on $M.$
\end{Theorem}

\begin{Example}
The natural para-hyperK\" ahler structure on the para-quaternionic vector space $\PH ^{n+1}$ is described in Example \ref{phn}.
Let the group $G = S^1$ of unit complex numbers  act on $\PH ^{n+1}$ by left multiplication, i.e.
$$
e^{it}\cdot h = e^{it}\cdot (h_1, \dots , h_{n+1}) :=(e^{it} h_1, \dots , e^{it} h_{n+1}), \enskip t \in \RR ,
$$
for $h =(h_1, \dots , h_{n+1}) \in \PH ^{n+1}.$ The action is isometric and preserves the para-hyperK\" ahler structure on $\PH ^{n+1}.$
We identify $\PH ^{n+1}$ with $\CC ^{2(n+1)} = \CC ^{n+1}\oplus \CC ^{n+1}$ by
$h =(z, w).$
After a straightforward calculation one finds that the equivariant  moment map of $G$ is
$$
f(z, w) = (|z|^2 + |w|^2 , -\Re (zw), -\Im (zw)).
$$
According to the notation of  Theorem \ref{glavna1}   choose  $\xi = (-1, 0, 0). $ Then
$$
{\cal K}_\xi  = \{ (z,w) \in \CC ^{2(n+1)} \mid  |z|^2 + |w|^2 = 1, zw =0 \} \subset \CC ^{2(n+1)}.
$$
The action of the group $G$ on ${\cal K}_\xi$ is given by  $e^{it}\cdot (z,w) = (e^{it}z, e^{-it}w), \enskip t \in \RR .$
After simple transformations one obtains that the resulting para-hyperK\" ahler manifold is a (real) submanifold of $\CC P^{2n + 1}$ of real codimension two, given in homogenous coordinates by
$$
M_\xi = \{ [z,w] \in \CC P^{2n + 1} \mid |z|^2 = |w|^2,   \Im (\bar zw) =0 \}.
$$
\end{Example}
\begin{Remark}
It is interesting that the resulting manifold is compact.
Analogous action  in the quaternionic case results in a  hyperK\" ahler structure on  the (noncompact) cotangent space $T^*\CC P^n$  of the complex projective space $\CC P^n .$
\end{Remark}

\subsection{Para-quaternionic K\" ahler reduction}

Let $M = ( M^{4n}, g, {\cal G})$ be a pa\-ra\--qua\-ternionic K\"
ahler manifold. Let a  Lie group  $G$ act  freely and
isometrically on $M$  and preserve the  $4$-form $\Omega = \Omega
({\cal G}),$ i.e. $g^*\Omega  = \Omega$ holds for a  $g \in G.$

\noindent
Denote by $V$ the unique Killing vector field corresponding to a Lie algebra vector $V^*\in {\cal g} = Lie (G).$
The section $\Theta _V$ of the  bundle $\Omega ^1({\cal G})$  of one forms with values in ${\cal G}$ defined by
$$
\Theta _{V^*}(X) := \sum _\alpha \omega _\alpha (V, X)J_\alpha,
$$
\noindent for $\enskip X \in TM,$ is well defined globally. Here
is the main reduction theorem.
\begin{Theorem}[Para-quaternionic K\" ahler reduction]
\label{PQKred} If a Lie group acts freely and isometrically on the
para-quaternionic K\" ahler manifold $M = (M^{4n}, g, {\cal G})$
and preserves the  $4$-form $\Omega (\cal G)$ then there exist a
unique section $f$ of bundle ${\cal g}^*\otimes {\cal G}$ such
that
\begin{equation}
\label{eq1}
\nabla f_{V^*} = \Theta _{V^*},
\end{equation}
for every $V^* \in {\cal g}.$ Moreover, the  group $G$ acts by isometries on the preimage ${\cal K} := f^{-1}(0)$ of the zero-section $0 \in {\cal g}^*\otimes {\cal G}.$ Suppose that  ${\cal K}$  is a smooth submanifold of $M$ and that the  quotient $\tilde M := {\cal K }/G$ has a smooth structure for which the projection $\pi: {\cal K} \to \tilde M$  is a pseudo-Riemannian submersion with induced metric $\tilde g$  on $\tilde M.$
Then the pseudo-Riemannian manifold $(\tilde M, \tilde g )$ is a para-quaternionic K\" ahler manifold with respect to the structure ${\cal G}'$ induced on $\tilde M$ from the structure ${\cal G}$ by the projection $\pi .$
\end{Theorem}
\begin{Remark}
Theorems \ref{glavna1} and \ref{PQKred} can be applied in case of a locally free action of the group $G.$ In that case the resulting para-quaternionic K\" ahler space may have an orbifold (rather than a manifold) structure (see \cite{GL}).
\end{Remark}
\noindent
For the proof of the Theorem \ref{PQKred} we need several lemmas.
\begin{Lemma}
\label{moment}
The unique solution $f_{V^*}, \enskip V^* \in {\cal g}$ of equation (\ref{eq1}) is given by
$$
f_{V^*} = \frac{4n}{K}\sum_\alpha \Tr (J_\alpha L_V)J_\alpha ,
$$
where $K$ is the scalar curvature of $M$ and $L_V = \nabla _V - {\cal L}_V$ is the Nomitzu operator (${\cal L}$ is the Lie derivative).
\end{Lemma}
\noindent {\bf Proof:} Notice that the scalar curvature of a
para-quaternionic K\" ahler manifold  is different from zero.
Since $V$ is a Killing vector field preserving para-quaternionic
K\" ahler structure  both operators ${\cal L}_V$ and  $\nabla _V,$
and hence the Nomitzu operator $L_V$ take values in the holonomy
$sp(n, \RR ) \oplus sp(1, \RR ).$ Let $e_1, \dots , e_{4n}$ be a
local pseudo-orthonormal reper. Then the relation
\begin{eqnarray*}
\nabla _k \Tr (J_\alpha L_V) & = & \nabla _k \biggl( \sum _{i,j = 1}^{4n}(J_\alpha)^i_j(V_{;i}^j)\biggr) = \\
& = & \sum _{i,j = 1}^{4n}(\nabla _ k J_\alpha)^i_jV_{;i}^j + \sum _{i,j = 1}^{4n} (J _\alpha)^i_j(\nabla _k V_{;i}^j) = \\
& = & \e _\alpha \sum _{i,j = 1}^{4n} (\mu _\gamma (e_k)J_\beta - \mu _\beta (e_k)J_\gamma)^i_jV_{;i}^j + \sum _{i,j = 1}^{4n} (J_\alpha)^i_j g(R(V, e_k)e_i, e_j)
\end{eqnarray*}
holds for any $k=1, \dots , 4n$ and  any $\alpha = 1, 2, 3.$  In other words, the relation
\begin{equation}
\label{cudo}
\nabla \Tr (J_\alpha L_V) = \e_\alpha (\mu_\gamma \Tr (J_\beta L_V) - \mu_\beta \Tr (J_\gamma L_V) ) + \Tr  (J_\alpha R(V, \cdot )),
\end{equation}
holds, where $(\alpha, \beta, \gamma)$ is a cyclic permutation of $(1, 2, 3).$
Summing up over $\alpha = 1, 2, 3$ most of the terms cancel
and one   obtains
\begin{equation}
\label{medju}
\nabla f_V = \frac{4n}{K}\sum \e_\alpha \Tr  (J_\alpha R(V, \cdot))J_\alpha.
\end{equation}
\noindent Since the curvature tensor $R$ of the manifold $M$ is of
type $sp(n, \RR ) \oplus sp(1, \RR )$, according to Theorem
\ref{decompSP} ii) the only component of the curvature tensor
which does not commute with all $J_\alpha$ is its projection  $R'$
onto ${\cal G}$
$$
R'(X, Y)  = -\frac{K}{(4n)^2}\sum _\alpha \e _\alpha g(J_\alpha X, Y)J_\alpha .
$$
Substituting $R'$ from the last relation into  relation (\ref{medju}) one easily checks
(\ref{eq1}). \qed

\begin{Lemma}
\label{indstruc}
The quotient space $\tilde M = ({\cal K}/G, \tilde g)$ from Theorem \ref{PQKred} has an almost  hermitian  para-quaternionic structure.
\end{Lemma}
\noindent
{\bf Proof :}
Let $V$ be a Killing vector field corresponding to a vector  $V^* \in {\cal g} = Lie (G).$
If  $e_k$ is a  vector tangent to ${\cal K }= f^{-1}(0)$ then using relation (\ref{cudo}) one obtains
$$
0 =\nabla _k f_{V^*} = \sum _\alpha \e _\alpha tr(J_\alpha  R(V, e_k))J_\alpha
 = \e _\alpha \frac{K}{4n} g(J_\alpha V, e_k)J_\alpha ,
$$
and hence
\begin{equation*}
\mbox{the vectors }\enskip J_\alpha V , \enskip \alpha = 1, 2, 3, \mbox{ are orthogonal to } {\cal K}.
\end{equation*}
Denote by $g' \subset TM$ the distribution spanned by all vector fields $V$ generated by the Lie algebra ${\cal g}.$
Then the subbundle ${\cal V} := \RR \langle g', J_1g', J_2g', J_3g' \rangle$ of $TM$  is ${\cal G}$ invariant. Moreover we have an orthogonal, ${\cal G}$ invariant splitting of $TM$
\begin{equation}
\label{splitting}
TM = {\cal V} \oplus {\cal V}^\perp
\end{equation}
which allows us to descend the para-quaternionic structure ${\cal G}$ of $M$ to  an almost para-quaternionic structure $\tilde{{\cal G}}$
of $ \tilde M,$ hermitian with respect to $\tilde g.$
We use the notation $\tilde J_\alpha$  and $\tilde \omega _\alpha ,\enskip  \alpha = 1, 2, 3,$ for  the
 basis of $\tilde{{\cal G}}$ and the corresponding $2$-forms.
Denote by $\tilde \Omega$ the $4$-form associated to $\tilde {\cal G}.$
\qed

\noindent In dimension four para-hyperK\" ahler manifolds are
characterized as pointwise Osserman manifolds which we are going
to define. Hence  we need the definition of  a (pointwise)
Osserman manifold. Let $(M, g)$ be a Riemannian or
pseudo-Riemannian manifold of any signature, and $R$ its curvature
tensor. We define its {\em Jacobi operator} in the unit direction
$X\in T_pM$ (i.e. $|X|^2 = \pm1$) by
$$
K_X(Y) := R(X, Y)X.
$$
It is a self-adjoint operator with respect to the metric $g.$ We
say that the manifold $M$ is {\em pointwise Osserman} if the
Jordan form of $K_X$ is independent of the  direction $X.$ We say
that $M$ is {\em (globally) Osserman} if  the Jordan form of $K_X$
is independent of both, the  direction $X\in T_pM$ and the point
$p$ of $M.$

\begin{Theorem}(\cite{ABB})
\label{oss1}
A pseudo-Riemannian manifold $(M, g)$ of signature $(2,2)$ is pointwise Osserman if and only if it is Einstein and self-dual.
\end{Theorem}

\begin{Lemma}
\label{oss2}
The $4$-dimensional manifold  $\tilde M = {\cal K}/G$ obtained by the para-quaternionic K\" ahler reduction (Theorem \ref{PQKred}) from a para-quaternionic K\" ahler manifold $M$ is a pointwise Osserman manifold.
\end{Lemma}
\noindent
{\bf Proof:}
We explain just the main steps. To simplify the notation we suppose that the group $G$ is one-dimensional and that $V$ is the Killing vector field of its action.
Denote by $V_X = \nabla _XV$ the covariant derivative of the vector field $V,$ where $\nabla$ is the Levi-Civita connection of the metric $g$ on $M.$
Let $i: {\cal K} \to M$ denote immersion, and let $\pi : {\cal K} \to \tilde M = {\cal K}/G$ be the natural projection.

\noindent
At first, using fundamental relations for immersion and submersion one  relates the curvature tensors
$R$ and $\tilde R$ of $M$ and $\tilde M,$ respectively.
Then, having a curvature tensor $\tilde R$ of $\tilde M,$ one can prove that the  Jacobi operator $\tilde K _{\tilde X}$ of the manifold $\tilde M$ in a point $\pi (u), \enskip u \in {\cal K},$
and in the unit direction $\tilde X \in T_{\pi (u)}\tilde M,$ is given by
\begin{eqnarray*}
\pi ^*(\tilde K_{\tilde X}\tilde Y)&  =&  \pi^*(\tilde R(\tilde X, \tilde Y )\tilde X) = h(R(X, Y) X) + \\
& + &
\frac{1}{|V|^2} \sum_{\alpha =1}^3 \e _\alpha
\left(g( X, J_\alpha (V_X))g( Y, h(J_\alpha (V_Y)))- g(Y, h(J_\alpha (V_X))) ^2\right)  + \\
&  + & \frac{1}{|V|^2} 3 g(Y,h(V_X)).
\end{eqnarray*}
Here $X=\pi ^* (\tilde X), \enskip Y=\pi ^* (\tilde Y)$ are horizontal lifts to $T_u{\cal K}$ of unit vectors  $\tilde X, \tilde Y \in T_{\pi (u)}M$  and $h(Z)$ denotes the ${\cal V}^\perp$ part (see (\ref{splitting})) of the vector $Z \in T_uM.$

\noindent
Using the decomposition of the curvature tensor $R$ from Theorem \ref{decompSP} ii) one obtains that
$$
h(R(X, Y), X) = \frac{K}{4n}Y,
$$
where $K$ is the  scalar curvature of $M.$ Since the Jacobi
operator $\tilde K _{\tilde X}$ is self-adjoint and the vector
$\tilde X$ is an eigenvector with eigenvalue $0$, we are
interested only in the restriction of $\tilde K _{\tilde X}$ to
the orthogonal complement of $\tilde X.$ One can calculate the
matrix of the restriction of $\tilde {\cal K} _{\tilde X}$  (for
example in the  basis $\tilde J_1 \tilde X, \tilde J_2 \tilde X,
\tilde J_3 \tilde X$) and show that it is diagonalizable with
eigenvalues
\begin{equation}
\label{eigen}
\lambda _1 = \lambda _2 = \frac{K}{4n} - \frac{2|h(V_X)|^2}{|V|^2}, \enskip
\lambda _2 = \frac{K}{4n} + \frac{4|h(V_X)|^2}{|V|^2} .
\end{equation}
One can show that the ratio $\frac{|h(V_X)|^2}{|V|^2}$ is independent of the  unit direction $X=\pi ^*(\tilde X)$ and hence the manifold $\tilde M$ is pointwise Osserman.
\qed

\begin{Remark}
In general the ratio $\frac{|h(V_X)|^2}{|V|^2}$ may depend of the  point $\pi (u), \enskip u \in {\cal K},$ and hence $\tilde M$ is not necessarily Osserman as we will see in Example \ref{nasprimer}.
\end{Remark}

\noindent
{\bf Proof of  Theorem \ref{PQKred}:} The existence of the unique   moment map $f \in {\cal g}^* \otimes \Omega ^0({\cal G})$ is proved in Lemma \ref{moment}.
\noindent
The moment map $f$ enjoys an important equivariance property in the following sense.
For any $g \in G,$ $V \in g$ and any point $p \in M$
$$
f_{g_*(V)}(p) = \tilde g(f_V(g^{-1}(p)))
$$
holds, where $\tilde g$ is the  map induced on ${\cal G}$ by $g.$
From this equivariance property it follows that the submanifold
$$
{\cal K} := \{ p \in M \mid f(p)  = 0 \} \subset M
$$
of $ M,$ is invariant under the (isometric) action of $G.$

\noindent
According to Lemma \ref{indstruc} there exists the natural  almost hermitian para-quaternionic  structure ${\cal G}$ on $\tilde M$ with $4$-form $\tilde \Omega .$
It remains to prove that $\tilde M$ is a para-quaternionic K\" ahler manifold.
If $dim \tilde M >4$ its enough to prove that $\tilde \Omega $ is parallel with respect to the Levi-Civita connection of the metric $\tilde g.$ The proof is the same as in the Riemannian case (see \cite{GL}).
If $dim \tilde M  = 4$ the proof follows from Lemma \ref{oss1} and Lemma \ref{oss2}. \qed
\begin{Remark}
In practice ${\cal K}$ must not be a differentiable submanifold of $M.$ In that case we can take a subset which is a submanifold of $M$ and which is invariant under the action of $G$ as shown in Example \ref{nasprimer}.
\end{Remark}
\begin{Example}
\label{nasprimer}
This example is from \cite{BV1} and it is a para-quaternionic version of the example by Galicki and Lawson (see \cite{GL}).
Let $p,q \in \NN$ be distinct, relatively prime natural numbers. We define the action of the Lie group $G := \{ e^{jt} \mid t \in \RR  \} \cong (\RR + )$ on $\PH P^2$   by
$$
\phi _{p,q}(t) \cdot [u_0, u_1, u_2] := [e^{jqt}u_0, e^{jpt}u_1
, e^{jpt}u_2],
$$
where $e^{jt} := \cosh t + j \sinh t$ and $[u_0, u_1, u_2]$  are homogenous coordinates on $\PH P^2.$ The action is free, isometric and preserves the para-quaternionic structure on $\PH P^2.$ One can show that the preimage of $0 \in \Imag \PH$ by the moment map
$f_{p,q} : \PH P^2 \to \Imag \PH$ is given by
$$
{\cal K}_{p,q}^0 = \{ [u_0, u_1, u_2] \in \PH P^2 \mid q\bar u_0ju_0 + p\bar u_1ju_1 + p\bar u_2ju_2 = 0 \}.
$$
The set ${\cal K}_{p,q}$ of regular points of ${\cal K}_{p,q}^0$ is
$$
{\cal K}_{p,q}= \{ [u_0, u_1, u_2] \in {\cal K}_{p,q}^0 \mid q^2 | u_0|^2 + p^2| u_1|^2+ p^2| u_2|^2 \neq  0 \}.
$$
The group $G$ acts freely and isometrically on ${\cal K}_{p,q}.$ The manifolds $\tilde M_{p,q} = {\cal K}_{p,q}/G$ obtained by the para-quaternionic reduction are $4$-dimensional, Einstein and self-dual  manifolds  of signature $(2,2).$ In \cite{BV2} it is proved by  computation of the eigenvalues (\ref{eigen}) that the manifolds $\tilde M_{p,q}$ are neither globally Osserman nor locally homogenous.
\end{Example}
\begin{Remark}
It is interesting that in case of action by $e^{it}$ instead of $e^{jt}$ one obtains an empty set as a preimage of $0 \in \Imag \PH $  by the moment map.
\end{Remark}

\section*{Acknowledgments}
\noindent I would like to express my gratitude to professor Dmitri
Alekseevsky for an invaluable help on the all stages of the work.
Special thanks to dr. Sandra Breimesser who read the first draft
of the manuscript and suggested many language improvements.



\begin{thebibliography}{10}

\bibitem{ABB}  D. Alekseevski, N. Bokan, N. Blazi\' c, Z. Raki\' c, {\it Self duality and pointwise Osserman condition}, Arch. mathematicum, 3(35), 193-201, (1999)

\bibitem{ABCV}  D.V. Alekseevski, N. Bla\v zi\' c, V.\ Cort\'es, S. Vukmirovi\'  c {\it A class of Osserman spaces }, in preparation

\bibitem{AC}  D.V. Alekseevski and V.\ Cort\'es, {\it
Classification of indefinite hyper-K\"ahler symmetric spaces}, to
appear in Asian J.\ Math., preprint 2000-79,
      Max-Planck-Institut f\"ur Mathematik, Bonn, math.DG/0007189.

\bibitem{AM}  D.V. Alekseevski, S. Marchifava {\em Quaternionic Structures on a Manifold and Subordinated Structures}, Anali di Mathematica pura ed applicata, (IV), Vol. CLXXI, 205-273, (1996)

\bibitem{AS} A. Swann, {\em HyperK\" ahler and quaternionic  K\" ahler geometry}, Math. Ann. {\bf 289}, 420-450 (1991)

\bibitem{B}  N. Blazi\' c, {\em Para-quaternionic projective spaces and pseudo-Riemannian geometry}, Publ. Inst. Mathem., 60(74), 101-107, (1996)

\bibitem{BBR} N.\ Bla\v zi\'c, N.\ Bokan,
Z.\ Raki\'c, {\it Osserman pseudo-Riemannian manifolds of
signature $(2,2)$}, to appear in Bull.\ Australian Math.\ Soc.\

\bibitem{BF} B. Feix, {\em HyperK\" ahler metrics on cotangent bundles}, J. reine angew. Math. (to appear)

\bibitem{BV1}  N. Blazi\' c, S. Vukmirovic, {\em Solutions of Yang-Mills equations on generalized Hopf bundles}, Journal of Geometry and Physics, 2001 (to appear)

\bibitem{BV2}  N. Blazi\' c, S. Vukmirovic, {\em Examples of self-dual, Einstein metrics of $(2,2)$ signature}, preprint

\bibitem{H1} N.J. Hitchin, {\em Hypersymplectic quotients},  Acta  Academiae Scientarum Tauriensis, Supplemento al numero {\bf 124} , 169-180, (1990)


\bibitem{CFG}  V. Cruceanu, P. Fortuny, and P. M. Gadea,
{\em A survey on para-complex geometry},
Rocky mountain J. of math., 26(1),  83--115,  1996.

\bibitem{GL} K. Galicki, B. Lawson {\em Quaternionic reduction and quaternionic orbifolds} Math. Ann. {\bf 282}, 1-21 (1988)

\bibitem{GMV} E.  Garcia-Rio, Y. Matsushita, R. Vazquez-Lorenzo, {\em Paraquaternionic K\" ahler manifolds}, to appear in Rocky Jour.

\bibitem{K} D.Kaledin, {\em HyperK\" ahler metrics on total spaces of cotangent bundles} in D. Kaledin, M. Verbitsky, {\em HyperK\" ahler manifolds} Math. Phys. Series, {\bf 12}, International Press, Cambridge MA, (1999)

\bibitem{Kam1} H. Kamada, {\em Neutral hyperK\" ahler structures on primary Kodaira surfaces}, Tsukuba J. Math. {\bf 23}, 321-332,  (1999)

\bibitem{Kam2} H. Kamada, {\em Neutral hyperK\" ahler structures on complex tori}, preprint 1998

\bibitem{HKLR} N.J. Hitchin, A. Karlhede, U. Lindst\" orm, M. Ro\v cek, {\em HyperK\" ahlet metrics and supersymmetry} Commun. Math. Phys. {\bf 108}, 535-589 (1987)

\bibitem{R} B. Rosenfeld, {\bf Geometry of Lie Groups}, Kluwer
Academic Publishers, (1997)

\bibitem{Sal} S. Salamon, {\bf Riemannian geometry and holonomy groups}, Longman Press, (1989)



\end{thebibliography}
\end{document}